\documentclass[11pt]{amsart}
\usepackage[utf8]{inputenc}
\usepackage[T1]{fontenc}
\usepackage{amsmath,amssymb,amsthm}
\usepackage[margin=1.15in]{geometry}
\usepackage[hidelinks]{hyperref}

\newtheorem{theorem}{Theorem}[section]
\newtheorem{proposition}[theorem]{Proposition}
\newtheorem{lemma}[theorem]{Lemma}
\newtheorem{corollary}[theorem]{Corollary}
\theoremstyle{remark}
\newtheorem{remark}[theorem]{Remark}
\theoremstyle{definition}

\newcommand{\ord}{\operatorname{ord}}
\newcommand{\jac}[2]{\left(\frac{#1}{#2}\right)}
\newcommand{\R}{\mathcal{R}}
\newcommand{\Sset}{\mathcal{S}}

\begin{document}

\title[On the number of elite primes]{On the number of elite primes}

\author{Ali Taheri}
\address{Department of Mechanical Engineering, Sharif University of Technology,
Tehran, Iran}
\email{ali.taheri82@sharif.edu}

\author{Amir Saleh}
\address{Department of Civil Engineering, Sharif University of Technology,
Tehran, Iran}
\email{amir.saleh01@sharif.edu}


\begin{abstract}
A prime $p$ is elite if all sufficiently large Fermat numbers $F_n=2^{2^n}+1$ are
quadratic nonresidues modulo $p$. Let $E(x)$ be the number of elite primes up to $x$.
The best published bounds are $E(x)\ll x/(\log x)^{3/2}$ unconditionally and
$E(x)\ll x^{5/6}$ under the Generalized Riemann Hypothesis (GRH), both due to Just.
Just and Pollack remarked, without a written proof, that sorting primes by the exact
power of $2$ dividing $p-1$ should give $x/(\log x)^{2-\varepsilon}$ and
$x^{3/4+\varepsilon}$. We give complete proofs of slightly sharper statements:
\[
E(x)\le\bigl(16+o(1)\bigr)\frac{x\log\log x}{(\log x)^2}
\quad\text{unconditionally,}\qquad
E(x)\ll x^{3/4}(\log x)^{1/2}\quad\text{under GRH.}
\]
We also recall, with short proofs, known facts about the period of $(F_n\bmod p)$
that go back to Aigner: every elite prime $p$ satisfies the defining condition from the
index $v_2(\ord_p 2)$ onwards, and the period is even unless $p\in\{3,5\}$. Finally, we show under GRH that, when the moduli involved are at most
$x^{1/3}$, the set of primes bounded by the sieve argument has size
$\asymp x\log\log x/(\log x)^2$. This indicates that the remaining factor $\log\log x$
cannot be removed by this method alone. All theorems of the paper
have been formalized in the Lean~4 proof assistant, with the Brun--Titchmarsh inequality
and two standard consequences of GRH as explicit hypotheses.
\end{abstract}

\maketitle

\section{Introduction}\label{sec:intro}

For $n\ge 0$ let $F_n=2^{2^n}+1$ be the $n$-th Fermat number. Following Aigner~\cite{Aig},
an odd prime $p$ is called \emph{elite} if there is an $n_0=n_0(p)$ such that
\[
\jac{F_n}{p}=-1\qquad\text{for every } n\ge n_0 .
\]
The elite primes below $7\cdot 10^5$ are
\[
3,\ 5,\ 7,\ 41,\ 15361,\ 23041,\ 26881,\ 61441,\ 87041,\ 163841,\ 544001,\ 604801,
\]
and Chaumont and
M\"uller~\cite{CM} listed all elite primes up to $2.5\cdot10^{11}$. Let
\[
E(x)=\#\{p\le x:\ p\ \text{is elite}\}.
\]

K\v{r}\'{\i}\v{z}ek, Luca and Somer~\cite[Theorem~4]{KLS} stated that $E(x)=O(x/(\log x)^2)$.
Just~\cite{Just} observed that their proof uses the estimate
$\prod_{i=0}^{2t}F_i<2^{2^{t+1}}$, which is false since the product equals
$2^{2^{2t+1}}-1$. With the correct value the argument gives
\begin{equation}\label{eq:just}
E(x)=O\!\left(\frac{x}{(\log x)^{3/2}}\right),
\end{equation}
and Just also proved $E(x)=O(x^{5/6})$ under GRH for Dirichlet $L$-functions.
In Remark~3.3 of~\cite{JP}, Just and Pollack stated without proof that keeping track of
the exact value of $e_p=v_2(p-1)$ improves these to
$O_\varepsilon(x/(\log x)^{2-\varepsilon})$ and, under GRH, $O_\varepsilon(x^{3/4+\varepsilon})$.
They carried out an argument of this type for anti-elite integers~\cite[\S3]{JP}.

The purpose of this note is to give complete proofs of these two claims, in a slightly
sharper form, together with some structural information that explains where the method
stops. The idea of splitting according to the exact value of $e_p$ is the one
indicated in~\cite[Remark~3.3]{JP}; our contribution is a complete and explicit account.

\begin{theorem}\label{thm:A}
As $x\to\infty$,
\[
E(x)\le \bigl(16+o(1)\bigr)\,\frac{x\log\log x}{(\log x)^2}.
\]
More precisely, for every integer $T\ge 1$ and every real $x$ with
$\bigl(2^{T}\cdot 2^{2^{T}}\bigr)^2\le x$,
\begin{equation}\label{eq:explicit}
E(x)\cdot 2^{T}\log x\le 4(T+1)\,x .
\end{equation}
\end{theorem}

\begin{theorem}\label{thm:B}
Assume GRH for all Dirichlet $L$-functions. Then
\[
E(x)=O\!\left(x^{3/4}(\log x)^{1/2}\right).
\]
\end{theorem}

Theorem~\ref{thm:A} uses only the Brun--Titchmarsh inequality in the form of Montgomery
and Vaughan~\cite{MV}. Theorem~\ref{thm:B} uses in addition the GRH bound for sums of a
non-principal character over primes. Nothing is assumed about the squarefreeness,
primality or factorization of Fermat numbers.

Theorem~\ref{thm:A} misses the bound stated in~\cite{KLS} by exactly a factor
$\log\log x$. Section~\ref{sec:barrier} explains why. All conditions used in the proof
are conditions on $p$ modulo a number of size $x^{O(1)}$. We show, under GRH, that when
this modulus is at most $x^{1/3}$, the primes satisfying all these conditions really
number $\asymp x\log\log x/(\log x)^2$ (Theorem~\ref{thm:barrier}). The proof of
Theorem~\ref{thm:A} uses moduli up to about $x^{1/2}$. This changes the parameter $T$ by
at most $O(1)$, and if the conditions remain equidistributed in this range, the argument
is sharp for the information it uses (Remark~\ref{rem:range}). An
improvement would need information about the Jacobi symbols $\jac{p}{F_n}$ for Fermat
numbers larger than any fixed power of $x$.

Section~\ref{sec:structure} concerns the period of the sequence $(F_n\bmod p)$. Write
$\ord_p(2)=2^{s_p}d_p$ with $d_p$ odd. By a theorem of Aigner~\cite{Aig}, in the form
given by Kla\v{s}ka~\cite[Theorem~20]{Kla}, the sequence is periodic exactly from the
index $s_p$ on, with period $\ord_{d_p}(2)$. Hence an elite prime satisfies
$\jac{F_n}{p}=-1$ already for every $n\ge s_p$ (and $s_p\le e_p$, with $s_p<e_p$ when
$e_p\ge3$). It is also known that the period is even for elite primes
$p>10$~\cite[\S2]{CM}. We state these facts as Theorem~\ref{thm:structure}, with short
proofs, and add one further necessary condition for eliteness (the second assertion of
Theorem~\ref{thm:structure}(c)). We also record the elementary bound
$\#\{p\le x: d_p\le y\}\le 2(y+1)\sqrt x$ (Lemma~\ref{lem:oddpart}). It shows that a
structural statement about $d_p$ for elite $p$ would give the bound of~\cite{KLS}.

Finally, all theorems of this paper have been checked in the Lean~4 proof assistant,
with the Brun--Titchmarsh inequality and the two consequences of GRH that we use as
explicit hypotheses; see Section~\ref{sec:lean}.

\section{Preliminaries}\label{sec:prelim}

Throughout, $p$ denotes an odd prime and $e_p=v_2(p-1)\ge1$. We write
$p-1=2^{e_p}k_p$ with $k_p$ odd, and $f_p=\ord_{k_p}(2)$ (with $f_p=1$ if $k_p=1$).
For $n\ge0$ let $\chi_n(a)=\jac{a}{F_n}$ be the Jacobi symbol. For a finite set
$B\subset\mathbb Z_{\ge0}$ put
\[
D_B=\prod_{n\in B}F_n,\qquad
\chi_B=\prod_{n\in B}\chi_n=\jac{\cdot}{D_B},\qquad \chi_\emptyset\equiv1 ;
\]
the second equality is the multiplicativity of the Jacobi symbol in the denominator.

\begin{lemma}\label{lem:fermat}
\begin{enumerate}
\item[(F1)] The Fermat numbers are pairwise coprime, and
$\prod_{i=0}^{n-1}F_i=2^{2^n}-1$ for $n\ge0$.
\item[(F2)] If a prime $r$ divides $F_n$, then $\ord_r(2)=2^{n+1}$; hence $e_r\ge n+1$.
\item[(F3)] Let $B\neq\emptyset$, let $M$ be a multiple of $D_B$, and let $\psi$ be a
Dirichlet character modulo an integer $Q$ coprime to $M$. Then $a\mapsto\psi(a)\chi_B(a)$,
restricted to integers coprime to $QM$, is a non-principal character modulo $QM$. In
particular $\chi_B$ induces a non-principal character modulo $M$.
\end{enumerate}
\end{lemma}

\begin{proof}
(F1) The product formula follows by induction from
$(2^{2^n}-1)F_n=2^{2^{n+1}}-1$. If $m<n$ and a prime $r$ divides $F_m$ and $F_n$, then
$r$ divides $F_n-\prod_{i<n}F_i=2$, which is impossible.

(F2) We have $2^{2^n}\equiv-1$ and $2^{2^{n+1}}\equiv1\pmod r$.

(F3) We have $F_0=3$, and for $n\ge1$,
$(2^{2^{n-1}})^2<F_n<(2^{2^{n-1}}+1)^2$, so no $F_n$ is a square. A product of pairwise
coprime positive integers is a square only if each factor is, so by (F1) $D_B$ is not a
square. Hence some prime $r$ divides $D_B$ to an odd power $v$. By the Chinese remainder
theorem choose $a$ with $a$ a quadratic nonresidue modulo $r$, $a\equiv1$ modulo
$D_B/r^{v}$, and $a\equiv 1$ modulo $Q$ and modulo every prime dividing $M$ but not $D_B$.
Then $\gcd(a,QM)=1$ and $\psi\chi_B(a)=\jac{a}{r}^{v}=-1$.
\end{proof}

We use the Brun--Titchmarsh inequality in the form of Montgomery and
Vaughan~\cite[Theorem~2]{MV}: for integers $1\le q<x$ and $a$ with $\gcd(a,q)=1$,
\begin{equation}\label{eq:BT}
\pi(x;q,a)\le\frac{2x}{\varphi(q)\log(x/q)} .
\end{equation}

The following reduction is standard (see~\cite[\S2]{Just}); we include the short proof.

\begin{lemma}\label{lem:reduction}
Let $p$ be elite and $e=e_p$. Then for every $n\ge e$ we have $p\nmid F_n$ and
\[
\chi_n(p)=\jac{p}{F_n}=\jac{F_n}{p}=-1 .
\]
Moreover $p\equiv1+2^{e}\pmod{2^{e+1}}$.
\end{lemma}

\begin{proof}
If $p\mid F_n$ with $n\ge e$, then (F2) gives $e_p\ge n+1>e$, a contradiction.
For $m\ge e$ we have $2^{m+f_p}-2^m=2^m(2^{f_p}-1)$, which is divisible by
$2^{e}k_p=p-1$. By Fermat's little theorem $F_{m+f_p}\equiv F_m\pmod p$ for all
$m\ge e$. Given $n\ge e$, choose $j$ with $n+jf_p\ge n_0(p)$; then
$\jac{F_n}{p}=\jac{F_{n+jf_p}}{p}=-1$. Since $n\ge e\ge1$, we have $F_n\equiv1\pmod4$,
and Jacobi reciprocity gives $\jac{p}{F_n}=\jac{F_n}{p}$. Finally $p-1=2^{e}k_p$ with
$k_p$ odd, so $p-1\equiv2^{e}\pmod{2^{e+1}}$.
\end{proof}

\begin{lemma}\label{lem:classes}
Let $1\le e\le T$ be integers. Put $I_e=\{e,e+1,\dots,T-1\}$,
\[
P_e=D_{I_e}=\prod_{n=e}^{T-1}F_n,\qquad M_e=2^{e+1}P_e ,
\]
and let $\R_e$ be the set of reduced residue classes $a\bmod M_e$ with
$a\equiv1+2^{e}\pmod{2^{e+1}}$ and $\chi_n(a)=-1$ for all $n\in I_e$. Then
\[
|\R_e|=\frac{\varphi(M_e)}{2^{T}}\qquad\text{and}\qquad P_e<2^{2^{T}} .
\]
\end{lemma}

\begin{proof}
By the Chinese remainder theorem, the class $a\bmod M_e$ corresponds to the pair
$(a\bmod 2^{e+1},\ a\bmod P_e)$, and the condition modulo $2^{e+1}$ fixes the first
component. Each $\chi_n(a)$ with $n\in I_e$ depends only on $a\bmod P_e$. For $a$
coprime to $P_e$ every $\chi_n(a)$ is $\pm1$, so the indicator of the remaining
conditions is
\[
\prod_{n\in I_e}\frac{1-\chi_n(a)}{2}=2^{-(T-e)}\sum_{B\subseteq I_e}(-1)^{|B|}\chi_B(a).
\]
Sum over the reduced residues $a\bmod P_e$. By (F3) and orthogonality, every term with
$B\neq\emptyset$ vanishes, which leaves $2^{-(T-e)}\varphi(P_e)$. Since
$\varphi(M_e)=2^{e}\varphi(P_e)$, we get $|\R_e|=\varphi(M_e)/2^{T}$. Finally, (F1) gives
$P_e\le\prod_{n=0}^{T-1}F_n=2^{2^T}-1$.
\end{proof}

\begin{remark}\label{rem:route}
The Lean formalization (Section~\ref{sec:lean}) proves $|\R_e|\cdot2^{T}=\varphi(M_e)$ by
a different route that avoids characters. The group of units modulo $M_e$ acts by
multiplication on the $2^{T}$ possible ``signatures''
$(a\bmod 2^{e+1},\ \{n\in I_e:\chi_n(a)=-1\})$, and this action is transitive. The
moving elements are built by the Chinese remainder theorem from the values $\chi_n(3)=-1$
($n\ge1$).
\end{remark}

\section{The unconditional bound}\label{sec:uncond}

\begin{proposition}\label{prop:explicit}
Let $T\ge1$ be an integer and $x$ a real number with $\bigl(2^{T}\cdot2^{2^{T}}\bigr)^2\le x$.
Then
\begin{align}
\#\{p\le x\ \text{elite}:\ e_p=e\}&\le\frac{4x}{2^{T}\log x}\qquad(1\le e\le T-1),
\label{eq:block}\\
\#\{p\le x\ \text{elite}:\ e_p\ge T\}&\le\frac{8x}{2^{T}\log x}.\label{eq:bigblock}
\end{align}
Consequently $E(x)\cdot2^{T}\log x\le4(T+1)\,x$.
\end{proposition}

\begin{proof}
Let $1\le e\le T-1$ and let $p\le x$ be elite with $e_p=e$. By
Lemma~\ref{lem:reduction}, $p$ is coprime to $P_e$ (and odd),
$p\equiv1+2^{e}\pmod{2^{e+1}}$, and $\chi_n(p)=-1$ for all $n\in I_e$. Hence
$p\bmod M_e\in\R_e$. By Lemma~\ref{lem:classes},
\[
M_e=2^{e+1}P_e<2^{T}\cdot2^{2^{T}}\le x^{1/2},
\]
so $M_e<x$ and $\log(x/M_e)>\tfrac12\log x$. Applying~\eqref{eq:BT} to each class in
$\R_e$ and using $|\R_e|=\varphi(M_e)/2^{T}$,
\[
\#\{p\le x\ \text{elite}:\ e_p=e\}\le\sum_{a\in\R_e}\pi(x;M_e,a)
\le\frac{\varphi(M_e)}{2^{T}}\cdot\frac{2x}{\varphi(M_e)\log(x/M_e)}
<\frac{4x}{2^{T}\log x}.
\]

If $e_p\ge T$, then $p\equiv1\pmod{2^{T}}$. Since $2^{T}\le x^{1/2}$, we have
$\log(x/2^{T})\ge\tfrac12\log x$, and $\varphi(2^{T})=2^{T-1}$. So~\eqref{eq:BT} gives
\[
\#\{p\le x:\ p\equiv1\ (\mathrm{mod}\ 2^{T})\}\le\frac{2x}{2^{T-1}\cdot\frac12\log x}
=\frac{8x}{2^{T}\log x}.
\]

Every elite prime has $e_p\ge1$ (the prime $2$ is not elite). Adding the $T-1$ bounds
\eqref{eq:block} and the bound~\eqref{eq:bigblock} gives
$E(x)\le(4(T-1)+8)\,x/(2^{T}\log x)=4(T+1)\,x/(2^{T}\log x)$.
\end{proof}

\begin{proof}[Proof of Theorem~\ref{thm:A}]
The explicit inequality~\eqref{eq:explicit} is Proposition~\ref{prop:explicit}. For
$x\ge2^{6}$, let $T\ge1$ be the largest integer with
$\bigl(2^{T}\cdot 2^{2^{T}}\bigr)^2\le x$, that is, with $2^{T+1}+2T\le\log_2x$. Then
$2^{T}\le\tfrac12\log_2x$, so $T\le\log_2\log_2x$. By maximality,
$2^{T+2}+2T+2>\log_2x$, so
\[
2^{T}>\tfrac14\bigl(\log_2x-2T-2\bigr)=\bigl(1+o(1)\bigr)\frac{\log x}{4\log2}.
\]
Proposition~\ref{prop:explicit} now gives
\[
E(x)\le\frac{4(T+1)\,x}{2^{T}\log x}
\le\bigl(16\log2+o(1)\bigr)\frac{(\log_2\log_2x)\,x}{(\log x)^2}
=\bigl(16+o(1)\bigr)\frac{x\log\log x}{(\log x)^2},
\]
because $\log2\cdot\log_2\log_2x=\log\log x-\log\log2$.
\end{proof}

\begin{remark}\label{rem:loss}
For each single value of $e$, the bound~\eqref{eq:block} has the order $x/(\log x)^2$
stated in~\cite{KLS}. The extra factor comes only from adding up the $\asymp\log\log x$
values of $e$ with $2^{e}\le\log x$. Each of them contributes the same amount. The
primes with $e_p=e$ have density $2^{-e}$, and the usable characters are $\chi_n$ with
$e\le n<T$, whose moduli already multiply to about $2^{2^{T}}\approx x^{1/2}$. So every $e$
saves exactly $2^{-e}\cdot2^{-(T-e)}=2^{-T}\asymp1/\log x$. By contrast, the proof of the
GRH bound in~\cite[Theorem~2]{Just} splits only into $e_p>t$ and $e_p\le t$, and uses the
characters $\chi_n$, $t\le n\le t+T$, for every $e_p\le t$. That costs both $x/2^{t}$ and
$x/2^{T}$ while the moduli are as large as $2^{2^{t+T+1}}$.
\end{remark}

\section{The bound under GRH}\label{sec:grh}

The only analytic input beyond~\eqref{eq:BT} is the following standard consequence of
GRH. It is also the input of~\cite[Theorem~2]{Just}.

\begin{lemma}\label{lem:grh}
Assume GRH for all Dirichlet $L$-functions. There is an absolute constant $c_0\ge1$ such
that for every $q\ge1$, every non-principal character $\chi\bmod q$ and every $x\ge2$,
\[
\Bigl|\sum_{p\le x}\chi(p)\Bigr|\le c_0\sqrt x\,\log(qx).
\]
\end{lemma}

\begin{proof}[Sketch of proof]
This is (13.21) in~\cite[Theorem~13.7]{MVbook}, which is stated for every modulus $q$
and every $x\ge2$, with no restriction on the size of $q$ relative to $x$. Since we
apply the lemma with moduli much larger than $x$, we sketch the argument. If $\log q\ge\sqrt x$ the bound is trivial. Otherwise let $\chi^*\bmod q^*$, $1<q^*\mid q$,
be the primitive character inducing $\chi$. The truncated explicit
formula~\cite[Ch.~19]{Dav} with $T=t$, the bound $\ll\log(q^*(u+2))$ for the number of zeros
with $u\le|\gamma|\le u+1$~\cite[Ch.~16]{Dav}, and $|t^{\rho}|=\sqrt t$ give
$\psi(t,\chi^*)\ll\sqrt t\,\log t\,\log(q^*t)$ for $2\le t\le x$ (when
$\log q^*\ge\sqrt t$ this is trivial). The prime powers $p^k$, $k\ge2$, contribute
$O(\sqrt t)$. Partial summation then gives
$\sum_{p\le x}\chi^*(p)\ll\sqrt x\log(q^*x)$. Finally, $\chi$ and $\chi^*$ differ only at
the at most $\log q/\log2$ primes dividing $q$.
\end{proof}

\begin{proof}[Proof of Theorem~\ref{thm:B}]
Let $x$ be large. Put $Y=x^{1/4}(\log x)^{-1/2}$ and let $T$ be the integer with
$2^{T}\le Y<2^{T+1}$. Then
\begin{equation}\label{eq:Tgrh}
2^{T}>Y/2,\qquad2^{T}\le x^{1/4},\qquad T<\log x,
\end{equation}
and $x/(Y\log x)=\sqrt x\,Y=x^{3/4}(\log x)^{-1/2}$.

\emph{Primes with $e_p\ge T$.} They satisfy $p\equiv1\pmod{2^{T}}$, and
$\log(x/2^{T})\ge\frac34\log x$. By~\eqref{eq:BT} their number is at most
\[
\frac{2x}{2^{T-1}\cdot\frac34\log x}=\frac{16\,x}{3\cdot2^{T}\log x}
<\frac{32}{3}\,x^{3/4}(\log x)^{-1/2}.
\]

\emph{Primes with $e_p=e$, $1\le e\le T-1$.} With $I_e$ as in Lemma~\ref{lem:classes} put
\[
w_e(p)=\prod_{n\in I_e}\frac{1-\chi_n(p)}{2}.
\]
Each factor is $0$, $\tfrac12$ or $1$, so $0\le w_e(p)\le1$ for every odd prime $p$, and
$w_e(p)=1$ if $p$ is elite with $e_p=e$ (Lemma~\ref{lem:reduction}). Every such $p$
satisfies $p\equiv1+2^{e}\pmod{2^{e+1}}$. Put
\[
S_B:=\sum_{\substack{p\le x\\ p\equiv1+2^{e}\ (2^{e+1})}}\chi_B(p)\qquad(B\subseteq I_e).
\]
Expanding the product,
\begin{equation}\label{eq:expand}
N_e:=\#\{p\le x\ \text{elite}:\ e_p=e\}
\le\sum_{\substack{p\le x\\ p\equiv1+2^{e}\ (2^{e+1})}}w_e(p)
=2^{-(T-e)}\sum_{B\subseteq I_e}(-1)^{|B|}S_B .
\end{equation}

For $B=\emptyset$, since $2^{e+1}\le2^{T}\le x^{1/4}$, \eqref{eq:BT} gives
\[
S_\emptyset=\pi(x;2^{e+1},1+2^{e})\le\frac{2x}{2^{e}\log(x/2^{e+1})}
\le\frac{8\,x}{3\cdot2^{e}\log x}.
\]

For $B\neq\emptyset$, orthogonality of the $2^{e}$ characters $\psi\bmod2^{e+1}$ gives
\[
S_B=2^{-e}\sum_{\psi\bmod2^{e+1}}\overline\psi(1+2^{e})\sum_{p\le x}\psi\chi_B(p).
\]
By (F3) with $M=D_B$ and $Q=2^{e+1}$, each $\psi\chi_B$ is a non-principal character
modulo $q_B=2^{e+1}D_B$. By Lemma~\ref{lem:classes}, $D_B\le P_e<2^{2^{T}}$, and since
$e+1\le T\le2^{T}$,
\[
\log q_B<(e+1+2^{T})\log2\le2^{T+1}\log2<2Y .
\]
Lemma~\ref{lem:grh} gives $\bigl|\sum_{p\le x}\psi\chi_B(p)\bigr|\le c_0\sqrt x\,(2Y+\log x)$
for each $\psi$, hence $|S_B|\le c_0\sqrt x\,(2Y+\log x)$.

There are fewer than $2^{T-e}$ nonempty $B\subseteq I_e$. So~\eqref{eq:expand}
and~\eqref{eq:Tgrh} give
\[
N_e\le\frac{8\,x}{3\cdot2^{T}\log x}+c_0\sqrt x\,(2Y+\log x)
<\Bigl(\frac{16}{3}+2c_0\Bigr)x^{3/4}(\log x)^{-1/2}+c_0\sqrt x\,\log x .
\]

\emph{Adding up.} By~\eqref{eq:Tgrh} there are fewer than $\log x$ values of $e$ in the
second case, so
\[
E(x)\le\Bigl(\frac{16}{3}+2c_0\Bigr)x^{3/4}(\log x)^{1/2}+c_0\sqrt x\,(\log x)^2
+\frac{32}{3}\,x^{3/4}(\log x)^{-1/2}\ll x^{3/4}(\log x)^{1/2}. \qedhere
\]
\end{proof}

\begin{remark}
Under GRH, using $\chi_B$ costs about $\sqrt x\log(\text{modulus})\approx\sqrt x\,2^{\max B}$,
while each level $e$ saves the density $2^{-T}$. Balancing $x/2^{T}$ against
$\sqrt x\,2^{T}$ gives $x^{3/4}$ per level. GRH alone gives no cancellation between the
prime sums of different characters $\chi_B$, so the exponent $3/4$ is the natural limit
of this method.
\end{remark}

\section{The period of \texorpdfstring{$(F_n\bmod p)$}{(F\_n mod p)} and the odd part of
\texorpdfstring{$\ord_p(2)$}{ord\_p(2)}}\label{sec:structure}

For an odd prime $p$ write $\ord_p(2)=2^{s_p}d_p$ with $d_p$ odd, and put
$f'_p=\ord_{d_p}(2)$ (so $f'_p=1$ if $d_p=1$).

Most of the next theorem is known. Part~(b) is a theorem of Aigner~\cite{Aig}; in the
form given here it is~\cite[Theorem~20]{Kla}, where Aigner's original statement is
reproduced as Theorem~1. The first assertion of~(c) follows from~(b) at once. That the
period is even for every elite prime $p>10$ is stated in~\cite[\S2]{CM}, with reference
to~\cite{Aig}. We include short proofs to keep the paper self-contained, and because
these statements are part of the formalization (Section~\ref{sec:lean}). We have not
found the second assertion of~(c) in the literature.

\begin{theorem}\label{thm:structure}
Let $p$ be an odd prime, $s=s_p$, $d=d_p$ and $f'=f'_p$.
\begin{enumerate}
\item[(a)] $s\le e_p$, and $s\le e_p-1$ if $e_p\ge3$.
\item[(b)] The sequence $(F_n\bmod p)_{n\ge s}$ is purely periodic with exact period $f'$.
\item[(c)] If $p$ is elite, then $\jac{F_n}{p}=-1$ for every $n\ge s$. If moreover
$s\ge1$, then also $\jac{1-2^{2^{s-1}}}{p}=-1$.
\item[(d)] If $p$ is elite and $p\notin\{3,5\}$, then $d>1$ and $f'$ is even.
\end{enumerate}
\end{theorem}

\begin{proof}
(a) We have $2^{s}\mid\ord_p(2)\mid p-1$, so $s\le e_p$. If $e_p\ge3$ then
$p\equiv1\pmod8$, so $2$ is a quadratic residue modulo $p$ and
$2^{(p-1)/2}\equiv1\pmod p$. Hence $2^{s}\mid\ord_p(2)\mid(p-1)/2$, and $s\le e_p-1$.

(b) Let $n\ge s$ and $k\ge1$. Then $F_{n+k}\equiv F_n\pmod p$ if and only if
$2^{s}d\mid2^{n}(2^{k}-1)$. Since $n\ge s$ and $d$ is odd, this holds if and only if
$d\mid2^{k}-1$, that is, if and only if $f'\mid k$.

(c) By (b) and the definition of elite primes, $\jac{F_n}{p}=-1$ for all $n\ge s$. Now let
$s\ge1$ and put $h=2^{s-1}d$. Then $\ord_p(2)=2h$, so $2^{h}\equiv-1\pmod p$. Write
$2^{f'}-1=dq$. Since $2^{f'}-1$ and $d$ are odd, $q$ is odd. With $n=s-1+f'\ge s$,
\[
2^{n}=2^{s-1}\bigl(1+dq\bigr)=2^{s-1}+hq,\qquad\text{so}\qquad
2^{2^{n}}=2^{2^{s-1}}\bigl(2^{h}\bigr)^{q}\equiv-2^{2^{s-1}}\pmod p .
\]
Hence $F_n\equiv1-2^{2^{s-1}}\pmod p$, and $\jac{1-2^{2^{s-1}}}{p}=\jac{F_n}{p}=-1$.

(d) Suppose first that $d=1$. Then $s\ge1$ and $\ord_p(2)=2^{s}$, so $p$ divides
$2^{2^{s}}-1$ but not $2^{2^{s-1}}-1$. Hence $p\mid F_m$ with $m=s-1$. If $m\ge2$, then
$e_p\ge m+1\ge3$ by (F2), so $\jac2p=1$. For every $n>m$ we have $2^{2^{n}}\equiv1$, hence
$F_n\equiv2$ and $\jac{F_n}{p}=1$, so $p$ is not elite (cf.~\cite[Theorem~2.1]{Mul}). Thus $d=1$ forces
$p\mid F_0F_1=15$, that is, $p\in\{3,5\}$.

Now let $d>1$ and put $y=2^{2^{s}}\in\mathbb F_p$. Since $\ord_p(2)=2^{s}d>2^{s}$, we have
$y\neq1$. By (b), $y^{2^{f'}}=2^{2^{s+f'}}=2^{2^{s}}=y$. Since $1+y^{2^{j}}\equiv F_{s+j}$,
telescoping gives
\[
\prod_{j=0}^{f'-1}F_{s+j}\equiv\prod_{j=0}^{f'-1}\bigl(1+y^{2^{j}}\bigr)
=\frac{y^{2^{f'}}-1}{y-1}=1\pmod p .
\]
By (c), each factor on the left is a quadratic nonresidue, so $(-1)^{f'}=\jac1p=1$.
\end{proof}

\begin{remark}
Part (c) shows that an elite prime satisfies the defining condition on the window
$n\ge s_p$, which by (a) can be strictly longer than the window $n\ge e_p$ from
Lemma~\ref{lem:reduction}, together with one extra condition at index $s_p-1$. We expect
that using $s_p$ instead of $e_p$ in the proof of Theorem~\ref{thm:A} changes only the
constant. The reason is that the primes with $s_p=s$ should have density $\asymp2^{-s}$,
so each level would again save only $\asymp2^{-T}$. We have not made this precise.
\end{remark}

\begin{remark}\label{rem:data}
For the twelve elite primes below $7\cdot10^{5}$, a direct computation gives
\[
d_p\in\{1,3,5,15,17,63,85\},\qquad f'_p\in\{1,2,4,6,8\},
\]
with $d_p=1$ only for $p=3,5$, in accordance with part (d). Heuristically, a prime
with orbit length $f'_p$ is elite with ``probability'' about $2^{-f'_p}$, which suggests
that elite primes have very small $d_p$.
\end{remark}

The next lemma shows that a statement of this kind would suffice for the bound stated
in~\cite{KLS}.

\begin{lemma}\label{lem:oddpart}
For real $x,y\ge1$,
\[
\#\{p\le x:\ d_p\le y\}\le2(y+1)\sqrt x .
\]
\end{lemma}

\begin{proof}
Fix $D\ge1$. Every prime $p$ with $\ord_p(2)=D$ divides $2^{D}-1$ and satisfies
$p\equiv1\pmod D$, so $p\ge D+1$. Their product divides $2^{D}-1<(D+1)^{D}$, so there are
fewer than $D$ of them. Those with $p\le x$ have the form $1+kD$ with
$1\le k\le x/D$, so there are at most $x/D$ of them. Hence, for fixed odd $d$,
\[
\#\{p\le x:\ d_p=d\}\le\sum_{s\ge0}\min\Bigl(2^{s}d,\ \frac{x}{2^{s}d}\Bigr)<4\sqrt x,
\]
because the terms with $2^{s}d\le\sqrt x$ form a geometric series with largest term
at most $\sqrt x$, and so do the remaining ones. There are at most $(y+1)/2$ odd $d\le y$.
\end{proof}

\begin{corollary}\label{cor:oddpart}
Suppose that $d_p\le\sqrt p/(\log p)^2$ for all sufficiently large elite primes $p$.
Then $E(x)\ll x/(\log x)^2$.
\end{corollary}

\begin{proof}
Split the elite primes in $(\sqrt x,x]$ into dyadic ranges $(z/2,z]$ with $z=x/2^{j}$.
In each range $d_p\le\sqrt z/(\log(z/2))^2$, so by Lemma~\ref{lem:oddpart} the range
contains $\ll z/(\log z)^2\ll 2^{-j}x/(\log x)^2$ elite primes. Summing over $j$ and adding
the at most $\sqrt x$ primes below $\sqrt x$ gives the claim.
\end{proof}

We do not know how to prove any nontrivial upper bound for $d_p$ for elite primes. The
Weil bound controls sums of $\jac{1+y}{p}$ over the whole subgroup of $\mathbb F_p^{*}$ of
order $d_p$, but not over the much thinner orbit
$\{2^{2^{n}}:n\ge s_p\}$, which has only $f'_p$ elements.

\section{Why the method stops at \texorpdfstring{$x\log\log x/(\log x)^2$}{x log log x/(log x)\^{}2}}
\label{sec:barrier}

For an integer $T\ge2$ put
\[
M_T=2^{T}\prod_{n=1}^{T-1}F_n ,
\]
and let $\Sset_T(x)$ be the set of odd primes $p\le x$ with $\chi_n(p)=-1$ for every $n$
with $e_p\le n<T$ (there is no condition if $e_p\ge T$). By Lemma~\ref{lem:reduction},
every elite prime $p\le x$ lies in $\Sset_T(x)$. The proof of
Proposition~\ref{prop:explicit} uses no other property of elite primes, so it shows
\[
\#\Sset_T(x)\cdot2^{T}\log x\le4(T+1)\,x\qquad\text{whenever }\bigl(2^{T}\cdot2^{2^{T}}\bigr)^2\le x .
\]
The next theorem shows that, under GRH, this is sharp up to a constant factor as long as
$M_T\le x^{1/3}$. For the slightly larger $T$ used in the proof of Theorem~\ref{thm:A},
see Remark~\ref{rem:range}.

\begin{theorem}\label{thm:barrier}
Assume GRH for all Dirichlet $L$-functions. If $T\ge2$ and $M_T\le x^{1/3}$, then
\[
\#\Sset_T(x)=\frac{T+1}{2^{T}}\,\mathrm{li}(x)+O\bigl(x^{5/6}\log x\bigr),
\]
with an absolute implied constant. In particular, if $T=T(x)$ is the largest integer
with $M_T\le x^{1/3}$, then $2^{T}\asymp\log x$ and
\[
\#\Sset_T(x)\asymp\frac{x\log\log x}{(\log x)^2}.
\]
\end{theorem}

\begin{proof}
Let $P=\prod_{n=1}^{T-1}F_n$, so $M_T=2^{T}P$. For a reduced residue $a\bmod M_T$ let
$e(a)=v_2(a-1)$ if $a\not\equiv1\pmod{2^{T}}$, and $e(a)=T$ otherwise; this depends only
on $a\bmod2^{T}$. Call $a$ admissible if $\chi_n(a)=-1$ for all $e(a)\le n<T$, and let
$\mathcal A_T$ be the set of admissible classes. For an odd prime $p\nmid M_T$ we have
$p\in\Sset_T(x)$ if and only if $p\le x$ and $p\bmod M_T\in\mathcal A_T$.

We count $\mathcal A_T$ with the Chinese remainder theorem
$(\mathbb Z/M_T)^{*}\cong(\mathbb Z/2^{T})^{*}\times(\mathbb Z/P)^{*}$. For
$1\le e\le T-1$ there are $2^{T-e-1}$ classes $b\bmod2^{T}$ with $e(b)=e$, and exactly one
with $e(b)=T$. For fixed $b$ with $e(b)=e<T$, the argument of Lemma~\ref{lem:classes}
(using (F3) with $M=P$) shows that exactly $2^{-(T-e)}\varphi(P)$ reduced classes
$c\bmod P$ satisfy $\chi_n(c)=-1$ for $e\le n<T$. Hence
\[
|\mathcal A_T|=\varphi(P)\Bigl(\sum_{e=1}^{T-1}2^{T-e-1}2^{-(T-e)}+1\Bigr)
=\frac{T+1}{2}\,\varphi(P)=\frac{T+1}{2^{T}}\,\varphi(M_T),
\]
since $\varphi(M_T)=2^{T-1}\varphi(P)$.

Under GRH, $\pi(x;q,a)=\mathrm{li}(x)/\varphi(q)+O(\sqrt x\log x)$ uniformly for
$q\le x$ and $\gcd(a,q)=1$ (see~\cite[Corollary~13.8]{MVbook}). Summing over $\mathcal A_T$ and adding
the $O(\log M_T)$ primes dividing $M_T$,
\[
\#\Sset_T(x)=\frac{T+1}{2^{T}}\,\mathrm{li}(x)+O\bigl(\varphi(M_T)\sqrt x\log x\bigr)
=\frac{T+1}{2^{T}}\,\mathrm{li}(x)+O\bigl(x^{5/6}\log x\bigr).
\]

For the second statement, note that $M_T=2^{T}(2^{2^{T}}-1)/3\ge2^{2^{T}}$ for $T\ge2$.
So $M_T\le x^{1/3}$ gives $2^{T}\le\frac13\log_2x$. By maximality,
$x^{1/3}<M_{T+1}<2^{T+1}\cdot2^{2^{T+1}}$, so $2^{T+1}+T+1>\frac13\log_2x$. Hence
$2^{T}\asymp\log x$ and $T\asymp\log\log x$. Then $(T+1)2^{-T}\mathrm{li}(x)\asymp
x\log\log x/(\log x)^2$, which dominates the error term $x^{5/6}\log x$.
\end{proof}

\begin{remark}\label{rem:numerics}
As a numerical check (not used anywhere), among the odd primes up to $5\cdot10^{6}$ the
proportion lying in $\Sset_T(5\cdot10^6)$ is $0.31235$ for $T=4$ and $0.18762$ for $T=5$.
The corresponding densities are $(T+1)2^{-T}=0.3125$ and $0.1875$. (Here $M_T$ exceeds
$x^{1/3}$, so these values lie outside the range of Theorem~\ref{thm:barrier}.)
\end{remark}

\begin{remark}\label{rem:range}
The Brun--Titchmarsh inequality~\eqref{eq:BT} requires the modulus to be smaller than
$x$, and $\log M_T\asymp2^{T}$. So in any argument of the type of
Proposition~\ref{prop:explicit} the parameter $T$ exceeds the one in
Theorem~\ref{thm:barrier} by at most $O(1)$, and $(T+1)2^{-T}$ changes by at most a
bounded factor. Assuming that equidistribution of the conditions persists for such $T$,
the factor $\log\log x$ in Theorem~\ref{thm:A} is therefore the true size of the sets
this method bounds. Removing it would require information about $\jac{p}{F_n}$ for primes
$p\le x$ and Fermat numbers $F_n$ larger than any fixed power of $x$. Equivalently, it
would require information about sums $\sum_{p\le x}\chi_B(p)$ of real characters whose
conductors exceed every fixed power of $x$. We know of no unconditional method that
gives such information, either for a single character or on average over the sparse
family $\{\chi_B\}$.
\end{remark}

\section{Formal verification}\label{sec:lean}

All theorems of this paper have been formalized in the Lean~4 proof
assistant~\cite{Lean4} on top of the \texttt{mathlib} library~\cite{mathlib}. Three analytic
inputs are not in \texttt{mathlib} and enter as explicit hypotheses of the theorems that
use them:
\begin{itemize}
\item the Brun--Titchmarsh inequality~\eqref{eq:BT} (stated for integer $x$, which
suffices);
\item Lemma~\ref{lem:grh}, in exactly the form stated there;
\item the prime number theorem for arithmetic progressions with the GRH error term used in
the proof of Theorem~\ref{thm:barrier}, namely
$\bigl|\pi(x;q,a)-\mathrm{li}(x)/\varphi(q)\bigr|\le C\sqrt x\log x$ for $x\ge2$,
$1\le q\le x$ and $\gcd(a,q)=1$, with $\mathrm{li}(x)=\int_2^x dt/\log t$.
\end{itemize}
GRH itself is therefore not stated in Lean; what is assumed is the two standard
consequences of it that the proofs use. The formalization compiles with no \texttt{sorry}
and no additional axioms. We ran \texttt{\#print axioms} on each of its 110 theorems; each
depends only on (a subset of) the three standard axioms of Lean (\texttt{propext},
\texttt{Classical.choice}, \texttt{Quot.sound}). It contains:
\begin{itemize}
\item Lemma~\ref{lem:fermat} (F1), (F2) and (F3), and the value $\chi_n(3)=-1$ for $n\ge1$;
\item Lemma~\ref{lem:reduction} in full;
\item Lemma~\ref{lem:classes}, in the form $|\R_e|\cdot2^{T}=\varphi(M_e)$ for
$1\le e\le T$, proved by the route of Remark~\ref{rem:route};
\item Proposition~\ref{prop:explicit} and inequality~\eqref{eq:explicit} for every real
$x$ with $\bigl(2^{T}\cdot2^{2^{T}}\bigr)^2\le x$, and the asymptotic form of
Theorem~\ref{thm:A}: for every $\varepsilon>0$,
$E(x)\le(16+\varepsilon)\,x\log\log x/(\log x)^2$ for all large $x$;
\item Theorem~\ref{thm:B} at integers $x$, with the explicit bound
$E(x)\le(16/3+3c_0)\,x^{3/4}(\log x)^{1/2}$ for $x\ge2^{72}$, together with the bound
for $N_e$ in its proof;
\item Theorem~\ref{thm:structure} in full. For (b), Lean proves that for $n\ge s_p$ one has
$F_{n+k}\equiv F_n\pmod p$ if and only if $f'_p\mid k$;
\item Lemma~\ref{lem:oddpart}, for all real $x,y\ge0$, and Corollary~\ref{cor:oddpart} at
integers $x$, with the explicit bound $E(x)\le192\,x/(\log x)^2$ for all large $x$;
\item Theorem~\ref{thm:barrier}. The second statement is proved in the explicit form: for
all large $x$, if $T\ge2$ and $M_T\le x^{1/3}<M_{T+1}$, then $\log x/12<2^{T}\le2\log x$
and $x\log\log x/(16(\log x)^2)\le\#\Sset_T(x)\le109\,x\log\log x/(\log x)^2$. The count
of admissible classes is done separately for each value of $e_p$, using
Lemma~\ref{lem:classes}.
\end{itemize}
The following are not formalized: the proof of Lemma~\ref{lem:classes} by characters given
above (Lean proves the same statement by a different route); the deduction of the three
hypotheses above from known results (Brun--Titchmarsh) and from GRH; and the remarks,
including the numerical data in Remarks~\ref{rem:data} and~\ref{rem:numerics}.
The Lean source files are available at
\url{https://github.com/tali7965/Elite-prime-LEAN}. They build with Lean
\texttt{v4.36.0-rc1} and \texttt{mathlib} commit \texttt{957aaac}.

\section*{Acknowledgements}
We thank Dr.~Matthew Just and Professor Paul Pollack. Their paper~\cite{JP}, in
particular Remark~3.3 there, led us to think about and work on this problem.

\section*{AI Declaration}
The authors used the large language model (LLM) in this
work. The LLM wrote the text of the paper. The
LLM also wrote the Lean~4 formalization described in Section~\ref{sec:lean}, which the
authors used to check the proofs. All theorems of the paper are formalized there, with
the three analytic inputs listed in that section as explicit hypotheses. The LLM also
wrote the Python code used to compute the numerical data in Remarks~\ref{rem:data}
and~\ref{rem:numerics}. The authors checked by hand the arguments that are not
machine-checked, and they take full responsibility for the content of the paper.

\section*{Software declaration}
The formal verification in Section~\ref{sec:lean} uses the Lean~4 proof
assistant~\cite{Lean4}, version \texttt{v4.36.0-rc1}, and the \texttt{mathlib}
library~\cite{mathlib}, commit \texttt{957aaac}. The Lean source files are available at
\url{https://github.com/tali7965/Elite-prime-LEAN}. The numerical data in
Remarks~\ref{rem:data} and~\ref{rem:numerics} were computed with Python.

\end{document}